\documentclass[12pt]{amsart}
\usepackage[cp1251]{inputenc}
\usepackage[english,russian]{babel}
\usepackage{amsmath}
\usepackage{amssymb}
\usepackage{amsfonts}
\usepackage{graphicx}

\title[Discrete singular integrals]
 {Discrete singular integrals in a half-space}

\author{Alexander V. Vasilyev, Vladimir B. Vasilyev}

\thanks{
This work was completed when the second author was a DAAD stipendiat and hosted in Institute of Analysis and Algebra, Technical University of Braunschweig.}

\begin{document}
\renewcommand{\refname}{References}
\renewcommand{\proofname}{Proof.}
\renewcommand{\abstractname}{Abstract}

\begin{abstract}
We consider Calderon -- Zygmund singular integral
 in the discrete half-space $h{\bf Z}^m_{+}$, where ${\bf Z}^m$ is entire lattice ($h>0$) in ${\bf R}^m$, and prove
 that the discrete singular integral operator is invertible in $L_2(h{\bf Z}^m_{+}$) iff such is its continual analogue.
 The key point for this consideration takes solvability theory of so-called periodic Riemann boundary problem, which is constructed by authors.


\end{abstract}
\maketitle

\section{Introduction}

We consider simplest Calderon-Zygmund operators of convolution type \cite{MP}
$$
v.p.\int\limits_{{\bf R}^m}K(x-y)u(y)dy=\lim\limits_{\stackrel{\varepsilon\to 0}{\scriptstyle N\to{+\infty}}}\int\limits_{\varepsilon<|x-y|<N}K(x-y)u(y)dy,
$$
where the kernel $K(x)$ satisfies the following conditions:

1) $K(tx)=t^{-m}K(x), \forall x\neq 0, t>0$;

2) $\int\limits_{S^{m-1}}K(\theta)d\theta=0,  S^{m-1}$ is the unit sphere in  ${\bf R}^m$;

3) $K(x)$ is differentiable on ${\bf R}^m\setminus\{0\}$.

Let us consider a discrete operator generated by the Calderon-Zygmund kernel $K(x),$
and defined on functions $u_h(\tilde{x}),$
$\tilde{x}\in h{\bf Z}^m,$ where ${\bf Z}^m$ is entire lattice ($h>0$)
in ${\bf R}^m,$ and the corresponding equation
\begin{equation}\label{1}
au_h(\tilde{x})+\sum\limits_{\tilde{y}\in h{\bf Z}_{+}^m} K(\tilde{x}-\tilde{y})
u_h(\tilde{y})h^m=v_h(\tilde{x}),\quad
\tilde{x}\in h{\bf Z}_{+}^m,
\end{equation}
$a$ is certain constant, in the discrete half-space
$h{\bf Z}_{+}^m=\left\{\tilde{x}\in h{\bf Z}^m:~\tilde{x_m}>0\right\},$
$u_h,v_h\in L_2(h{\bf Z}_{+}^m).$

By definition we put $K(0)=0,$ and for the operator
$$
u_h(\tilde{x})\mapsto au(\tilde{x})+\sum\limits_{\tilde{y}\in h{\bf Z}^m}
 K(\tilde{x}-\tilde{y})
u_h(\tilde{y})h^m,\quad
\tilde{x}\in h{\bf Z}^m,
$$
we introduce its symbol by the formula
$$
\sigma_h(\xi)=a+\sum\limits_{\tilde{x}\in h{\bf Z}^m}e^{-i\xi\tilde{x}}K(\tilde{x})h^m;
$$
it is periodic function with basic cube period $[-\pi h^{-1};~\pi h^{-1}]^m.$

The sum for $\sigma_h(\xi)$ is defined as a limit of partial sums over cubes $Q_N$
$$
\lim\limits_{N\to \infty}\sum\limits_{\tilde{x}\in Q_N}
e^{-i\xi\tilde{x}}K(\tilde{x})h^m,
$$
$$
Q_N=\left\{\tilde{x}\in h{\bf Z}^m:~|\tilde{x}|\le  N,~|\tilde{x}|=
\max\limits_{1\le k \le m}|\tilde{x}_k|\right\}.
$$

It is very similar classical symbol of Calderon-Zygmund operator \cite{MP}, which is defined as Fourier transform
of the kernel $K(x)$
in principal value sense
$$
\sigma(\xi)=\lim\limits_{\stackrel{N\to\infty}{\scriptstyle \varepsilon \to 0}}
\int\limits_{\varepsilon <|x|<N}K(x)e^{i\xi x}dx.
$$

Key point of our study is a theorem proved in \cite{V1}, asserting that images of $\sigma$ and $\sigma_h$
are the same.

We also introduce continual equation in a half-space
\begin{equation}\label{2}
au(x)+\int\limits_{{\bf R}^m_+}K(x-y)u(y)dy=v(x),~x\in{\bf R}^m_+,
\end{equation}
and we'll prove that the equations \eqref{1} and \eqref{2} are uniquely solvable or unsolvable simultaneously for all $h>0$ in corresponding spaces.

\section{Discrete Calderon-Zygmund Operators}

\subsection{Symbol properties}

We recall some properties of symbols $\sigma(\xi)$ and $\sigma_h(\xi)$, which are needed for us \cite{V1}.

{\bf Lemma 1.}
{\it $\lim\limits_{h\to 0}\sigma_h(\xi)=\sigma(\xi),~~\forall\xi\neq 0$.}

\begin{proof}
 Indeed, if we fix $\xi\neq 0$, then by definition of integral as a limit of integral sums, we finish the proof.
\end{proof}

{\bf Lemma 2.}
{\it $\sigma_h(\xi)=\sigma_1(h\xi),~~\forall h>0,~~\xi\in [-\pi h^{-1},\pi h^{-1}]^m$.}

\begin{proof}
$$
\sigma_h(\xi)=\sum_{\tilde x\in h{\bf Z}^m}K(\tilde x)e^{-i\tilde x\cdot\xi}h^m=
$$
$$
=\sum\limits_{\tilde y\in {\bf Z}^m}K(h\tilde y)e^{-i\tilde y\cdot h\xi}h^m=\sum\limits_{\tilde y\in {\bf Z}^m}K(\tilde y)e^{-i\tilde y\cdot h\xi}=\sigma_1(h\xi).
$$
\end{proof}

{\bf Lemma 3.}
{\it  The images of $\sigma$ and $\sigma_h$ are the same, and their values are constant for any ray from origin.}

\begin{proof}
 It follows from previous lemmas immediately, because if we fix $\xi$, then $\sigma_1(0)=\sigma(\xi)\Longrightarrow\sigma_h(0)=\sigma(\xi).$
\end{proof}

\subsection{Symbols and operators}

We consider more general case in the whole space ${\bf R}^m$
$$
(M_1P_++M_2P_-)U=V,
$$
taking into account that $M_1,M_2$ are operators of type \eqref{2}, and $P_+,P_-$ are restriction operators on ${\bf R}^m_{\pm}=\{x=(x_1,...,x_m),~\pm x_m>0\}$.
It is easily verified that the equation \eqref{2} is a special case for such equation, when $M_2\equiv I, I$ is identity operator.

If we'll denote the Fourier transform by letter $F$, and use the notations \cite{GK}
$$
FP_+=Q_{\xi'}F,~FP_-=P_{\xi'}F,
$$
$$
P_{\xi'}=1/2(I+H_{\xi'}),~Q_{\xi'}=1/2(I-H_{\xi'}),
$$
where $H_{\xi'}$ is Hilbert transform on variable $\xi_m$ for fixed $\xi'=(\xi_1,...,\xi_{m-1})$ \cite{Es}:
$$
(H_{\xi'}u)(\xi',\xi_m)\equiv\frac{1}{\pi i}v.p.\int\limits_{-\infty}^{+\infty}\frac{u(\xi',\tau)}{\tau-\xi_m}d\tau,
$$
then the equation mentioned after applying the Fourier transform will be the following equation with the parameter $\xi'$:
$$
\frac{\sigma_{M_1}(\xi',\xi_m)+\sigma_{M_2}(\xi',\xi_m)}{2}\tilde U(\xi)+
$$
$$
+\frac{\sigma_{M_1}(\xi',\xi_m)-\sigma_{M_2}(\xi',\xi_m)}{2\pi i}v.p.
\int\limits_{-\infty}^{+\infty}\frac{\tilde U(\xi',\eta)}{\eta-\xi_m}d\eta=\tilde V(\xi)
$$
( $\tilde{}$  denotes the Fourier transform).

This equation is closely related to boundary Riemann problem with the parameter $\xi'$) with coefficient \cite{M,G}
$$
G(\xi',\xi_m)=\sigma_{M_1}(\xi',\xi_m)\sigma_{M_2}^{-1}(\xi',\xi_m).
$$

\section{Periodic Riemann Boundary Problem}

The theory of periodic Riemann boundary problem was constructed by authors \cite{V2} (see also forthcoming paper with the same name in {\it Differential Equations}) with full details, and now we will use
its general consequences.

Let's denote ${\bf Z}_+=0,1,2,..., {\bf Z}_-={\bf R}\setminus{\bf Z}_+$. The Fourier transform for function of discrete variable is the series
\begin{equation}\label{3}
(Fu)(\xi ) = \sum\limits_{k = - \infty }^{ + \infty } {u(k)e^{ - ik\xi }},~\xi\in[-\pi,\pi].
\end{equation}

Let's consider the Fourier transform \eqref{3} for the indicator of ${\bf Z}_+$:
 $$
 \chi _{{\bf Z}_ + } (x) = \left\{ {{\begin{array}{*{20}c}
 {1,} \hfill & {x \in {\bf Z}_ + } \hfill \\
 {0,} \hfill & {x \notin {\bf Z}_ + } \hfill \\
\end{array} }}. \right.
$$
For summable functions their product transforms to convolution of their Fourier images on the segment $[-\pi,\pi]$ but for our case
$F(\chi _{_{{\bf Z}_ + } } \cdot u)$ one of functions $\chi_{{\bf Z}_ +}$ is not summable.Thus, first we introduce some regularizing multiplier and evaluate the following Fourier transform

\[
F(e^{ - \tau \,k} \cdot \chi _{_{{\bf Z}_ + } } )(\xi)\mathop =
\frac{1}{2\pi }\sum\limits_{k \in {\bf Z}_ + } {e^{-\tau k}e^{ - ik\xi }} =
\frac{1}{2\pi }\sum\limits_{k \in {\bf Z}_ + } {e^{-\tau k - ik\xi }} =
\]

\[
 = \frac{1}{2\pi }\sum\limits_{k \in{\bf Z}_ + } {e^{ - ik(\xi + i\tau )}} =
\frac{1}{2\pi }\sum\limits_{k \in{\bf Z}_ + } {e^{ - ikz}}
,~
\quad
\tau \to 0,~
\quad
z = \xi + i\tau ,\quad \tau > 0.
\]

The Fourier transform for the function $u(n)$ we'll denote $\hat {u}(\xi )$, it is left to find the sum for
$e^{ - ikz}$,

\[
\frac{1}{2\pi }\sum\limits_{k \in{\bf Z}_ + } {e^{ - ikz}} = \frac{1}{2\pi }(1 +
e^{ - iz} + e^{ - 2iz} + ...) = \frac{1}{2\pi }\frac{1}{1 - e^{ - iz}},
\]

After some transformations:

\[
F(\chi_{{\bf Z}_+}\cdot u)(\xi)\mathop = \lim\limits_{\tau \to 0 +
} \left( {\frac{\hat {u}(\xi )}{4\pi } + \frac{1}{4\pi i}\int\limits_{ - \pi
}^\pi {\hat {u}(t)\cot\frac{z - t}{2}dt} } \right),\quad z = \xi + i\tau .
\]

According to Sokhotskii formulas (these are almost same for periodic kernel $\cot(x)$) (see also classical books \cite{M,G})

\[
F(\chi_{{\bf Z}_+}\cdot u)(\xi ) = \frac{\hat {u}(\xi )}{4\pi } + \frac{1}{4\pi
i}\int\limits_{ - \pi }^\pi {\hat {u}(t)\cot\frac{\xi - t}{2}dt} + \frac{\hat
{u}(\xi )}{2}
\]

If we introduce the function $\chi _{_{{\bf Z}_ - } } (x)$
and consider the Fourier transform for the product
$F(\chi _{{\bf Z}_ - } \cdot u)$ with preliminary regularization, then we have

\[
F(e^{ - \tau \,k} \cdot \chi _{_{Z_ - } } )\mathop =
\frac{1}{2\pi }\sum\limits^{ - 1}_{ - \infty } {e^{\tau k}e^{ - ik\xi }}
= \frac{1}{2\pi }\sum\limits^{ - 1}_{ - \infty } {e^{\tau k - ik\xi }} =
\]

\[
 = \frac{1}{2\pi }\sum\limits^{ - 1}_{ - \infty } {e^{ - ik(\xi + i\tau
)}} = \frac{1}{2\pi }\sum\limits^{ - 1}_{ - \infty } {e^{ - ikz}}
,~
\quad
\tau \to 0,
\quad
z = \xi + i\tau ,\quad \tau < 0.
\]

Further,

\[
\frac{1}{2\pi }\sum\limits^{ - 1}_{ - \infty } {e^{ - ikz}} =
\frac{1}{2\pi }( - 1 + 1 + e^{iz} + e^{2iz} + ...) = - \frac{1}{2\pi } +
\frac{1}{2\pi }\frac{1}{1 - e^{iz}},
\]

With the help of some elementary calculations:

\[
F(\chi_{{\bf Z}_-}\cdot u)\mathop = \limits_{n \to \xi } \mathop {\lim }\limits_{\tau \to 0}
\left( { - \frac{\hat {u}(\xi )}{4\pi } - \frac{1}{4\pi i}\int\limits_{ -
\pi }^\pi {\hat {u}(t)\cot\frac{z - t}{2}dt} } \right),\quad z = \xi + i\tau
.
\]

Applying Sokhotskii formulas, we have:

\[
F(\chi_{{\bf Z}_-}\cdot u)(\xi ) = \frac{\hat {u}(\xi )}{4\pi } + \frac{1}{4\pi
i}\int\limits_{ - \pi }^\pi {\hat {u}(t)\cot\frac{\xi - t}{2}dt} + \frac{\hat
{u}(\xi )}{2}.
\]

To verify one can find the sum for $F(\chi_{{\bf Z}_+}\cdot u)$, $F(\chi_{{\bf Z}_-}\cdot u)$ and obtain:

\[
F(\chi_{{\bf Z}_+}\cdot u)+F(\chi_{{\bf Z}_-}\cdot u) = \frac{\hat {u}(\xi )}{4\pi } +
\frac{1}{4\pi i}\int\limits_{ - \pi }^\pi {\hat {u}(t)\cot\frac{\xi -
t}{2}dt} + \frac{\hat {u}(\xi )}{2} -
\]

\[
 - \frac{\hat {u}(\xi )}{4\pi } - \frac{1}{4\pi i}\int\limits_{ - \pi }^\pi
{\hat {u}(t)\cot\frac{\xi - t}{2}dt} + \frac{\hat {u}(\xi )}{2}=\hat {u}(\xi ).
\]

These calculations lead to certain periodic Riemann boundary value problem, for which the solvability conditions are defined by the index of its coefficient. The problem is formulated
as following way: finding two functions $\Phi^{\pm}(t)$ which admit an analytical continuation into upper and lower half-strip in the complex plane $\bf C$, real part is the segment $[-\pi,\pi]$,
and their boundary values satisfy the relation
$$
\Phi^+(t)=G(t)\Phi^-(t)+g(t),~~~t\in[-\pi,\pi],
$$
$G(t), g(t)$ are given functions on $[-\pi,\pi]$, and such that $G(-\pi)=G(\pi), g(-\pi)=g(\pi)$.
Index for such problem is called the integer number
$$
\ae=\frac{1}{2\pi}\int\limits_{-\pi}^{\pi}d\arg G(t).
$$

\newpage
\centerline{\includegraphics*[width=15in, height=12in, keepaspectratio=false]{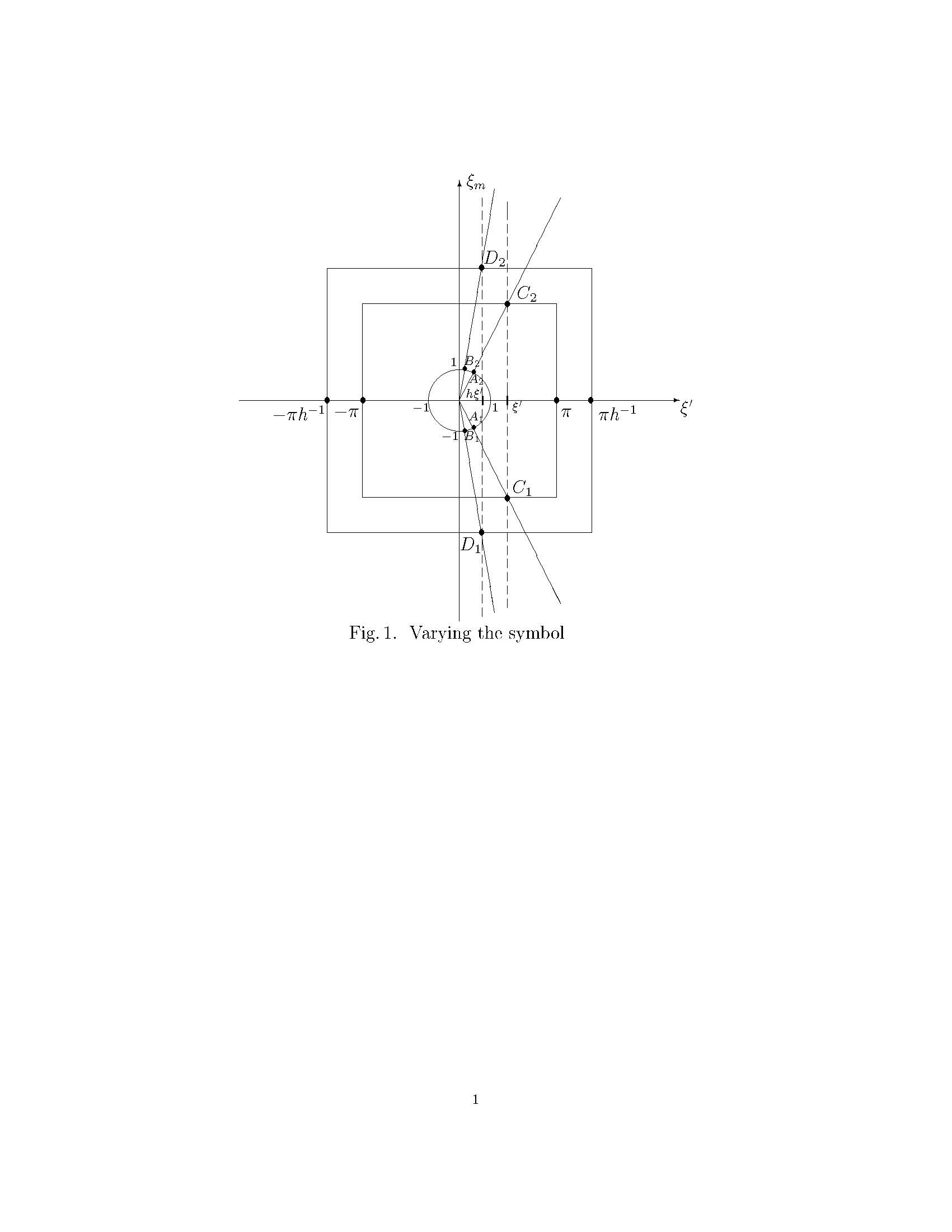}}


\section{Solvability Conditions}

Here we suppose additionally, that the symbol $\sigma(\xi',\xi_m)$ satisfies the condition
$$
\sigma(0,...,0,-1)=\sigma(0,...,0,+1).
$$

{\bf Main Theorem.}
{\it The equations \eqref{1} and \eqref{2} are uniquely solvable or unsolvable simultaneously for all $h>0$.}

\begin{proof}
 We need to look our symbols $\sigma(\xi)$ and $\sigma_h(\xi)$ more exactly. We'll illustrate our consideration with the help of Fig.1.

If we fix $\xi'$ in the cube $[-\pi,\pi]^m$, then under varying $\xi_m$ on $[-\pi,\pi]$ the argument of $\sigma_1(\xi)$ will vary along the curve on a
cubical surface of $[-\pi,\pi]^m$, which unites the points $C_1$ and $C_2$ (for the case $m\geq 3$ all such curves are homotopic, and for the case $m=2$ there are two curves left and right one).

This varying corresponds to the varying of the argument of function $\sigma(\xi)$ along the curve from point $A_1$ to point $A_2$ on the unit sphere. Further, if we consider the symbol $\sigma_h(\xi)$
now on the cube $[-h^{-1}\pi,h^{-1}\pi]^m$, then according to lemma 2.2 $h\xi_m$ will be varied on $[-\pi,\pi]$ also under fixed $h\xi'$. In other words, the argument of
$\sigma_h(\xi)$ for fixed $\xi'$  (we consider small $h>0$) will be varied along a curve on a cubical surface of $[-h^{-1}\pi,h^{-1}\pi]^m$, which unites the points $D_1$ and $D_2$. It corresponds to
varying argument of function $\sigma(\xi)$ from point $B_1$ to point $B_2$ on the unit sphere. Obviously, under decreasing $h$ the sequence $A_1,B_1,...$ will be convergent to the south pole of the unit sphere
$(0,...,0,-1)$, and the sequence $A_2,B_2,...$ to the north pole $(0,...,0,+1)$. Thus, because the variation of an argument of $\sigma_h(\xi)$ on $\xi_m$ under fixed $\xi'$ is $2\pi k$ ($\sigma_h$ is periodic function),
then under additional assumption
$$
\sigma(0,...,0,-1)=\sigma(0,...,0,+1)
$$
(this property is usually called the transmission property)
we'll obtain that variation of an argument for the function $\sigma(\xi)$ under varying $\xi_m$ from $-\infty$ to $+\infty$ under fixed $\xi'$ (this variation of $\sigma(\xi)$ moves along the arc of a big half-circumference on the unit sphere)
is also $2\pi k$. According to our assumptions on continuity of $\sigma(\xi)$ on the unit sphere, it will be the same number $2\pi k, k\in{\bf Z}$,
$$
\lim\limits_{h\to 0}\int\limits_{-\pi h^{-1}}^{\pi h^{-1}}d\arg\sigma_h(\xi',\xi_m)=\int\limits_{-\infty}^{+\infty}d\arg\sigma(\xi',\xi_m),~~~~~~\forall\xi'\neq 0.
$$

So, both for the equation \eqref{1} and equation \eqref{2} the uniquely solvability condition is defined by the same number. This completes the proof.
\end{proof}

\section{Conclusion}

We see that both continual and discrete equations are solvable or unsolvable simultaneously, and then we need to find good finite approximation for infinite system of linear algebraic equations for computer
calculations. First steps in this direction were done in the paper \cite{V3}, where the authors suggested to use fast Fourier transform.


\subsection*{Acknowledgment}
Many thanks to DAAD and Herr Prof. Dr. Volker Bach for their support.

\vspace{1cm}

Alexander V. Vasilyev\\
Department of Mathematical Analysis\\
Belgorod State University\\
Studencheskaya 14/1\\
Belgorod 308007, Russia\\
\vspace{1cm}
\email{alexvassel@gmail.com}\\

Vladimir B. Vasilyev\\
\address{Chair of Pure Mathematics\\
Lipetsk State Technical University\\
Moskovskaya 30\\
Lipetsk 398600, Russia\\
vbv57@inbox.ru

\end{document}